\documentclass[11pt]{article}
\usepackage{lscape}
 \usepackage[ansinew]{inputenc}
 \usepackage[dvips]{graphicx}
\usepackage[psamsfonts]{amssymb} \usepackage[all]{xy}  \usepackage[all]{xy}
 \usepackage[psamsfonts]{eucal}
 \usepackage[psamsfonts]{eucal}
\usepackage{amssymb, amsmath}
\usepackage{amsthm}
\usepackage{geometry}
\usepackage{enumerate}
\usepackage{float}
\newcommand{\be}{\begin{equation}}
\newcommand{\ee}{\end{equation}}

\newcommand{\bd}{\begin{displaymath}}
\newcommand{\ed}{\end{displaymath}}
\newcommand{\ben}{\begin{enumerate}}
\newcommand{\een}{\end{enumerate}}

\title{The depth of a Riemann surface and of a right-angled Artin group} 
\author{Yves F\'elix and Steve Halperin}
\date{\today}

\begin{document}

\maketitle

\begin{abstract}
We consider two families of spaces, $X$ : the  closed orientable Riemann surfaces of genus $g>0$ and  the classifying spaces of    right-angled Artin groups. In both cases we   compare   the depth of the fundamental group with the depth of  an associated  Lie algebra, $L$, that can be determined by the minimal Sullivan algebra. For these spaces we prove that 
$$  \mbox{depth} \,\mathbb Q[\pi_1(X)] = \mbox{depth}\, {L}\,$$
and give precise formulas for the depth.
    \end{abstract}


\section{Introduction}

 In this paper we show that the invariants \emph{depth}, defined in three distinct categories (Sullivan algebras, augmented and possibly graded algebras, $R$; and graded Lie algebras, $E$) all coincide in the case of compact oriented Riemann surfaces and right-angled Artin groups.
 
 \emph{Throughout the base field is $\mathbb Q$}, and we adopt the following notation: $UE $ is the universal enveloping algebra, and  
 $\widehat{R} := \varprojlim_n R/I^n $ is the completion of $R$
 (Here $I^n$ is the $n^{th}$ power of the augmentation ideal.)
 
 Then depth$\, R$ and depth$\, E$ are defined by 
 $$\mbox{depth$\,R=$ least $p$ (or $\infty$) such that Ext}_R^p(\mathbb Q, R)\neq 0$$
 and
 $$\mbox{depth$\, E=$ least $p$ (or $\infty$) such that Ext}_{UE}^p(\mathbb Q, \widehat{UE})\neq 0\,.$$

On the other hand, associated with any group $G$ is its \emph{Malcev Lie algebra}, (\cite{Mal}), 
$$M_G = \bigoplus_{n\geq 1} M_G(n)\,, \hspace{1cm} M_G(n) = G^n/G^{n+1} $$
where the normal subgroups $G^n$ are defined inductively: $G^{n+1}$ is generated by the commutators $[a,b]=aba^{-1}b^{-1}$ with $a\in G$ and $b\in G^n$. The Lie bracket in $M_G$ is induced by $[\,\,, \,\,]: G^p\times G^q \to G^{p+q}$. The   Lie algebra $$L_G:= M_G\otimes \mathbb Q$$
is the \emph{rational Malcev Lie algebra}.

More generally, a Lie algebra $E$ with a \emph{weight decomposition} (briefly, \emph{a weighted Lie algebra}) is a Lie algebra of the form $E = \oplus_{n\geq 1} E(n)$ satisfying $[E(p), E(q)]\subset E(p+q)$. Since by definition $[L_G(p), L_G(q)]\subset L_G(p+q)$ it follows that the rational Malcev Lie algebra is a Lie algebra with a weight decomposition.

Also associated with $X$ is its minimal Sullivan model $(\land V_X,d)$. The definition of a minimal Sullivan algebra $(\land V,d)$ and its depth, are recalled in section 2, and we note that by \cite[Theorem 10.1]{FHTII} 
$$\mbox{depth}\,(\land V_X,d)\leq \,\mbox{ cat}\, X,$$ where cat$\,X$ is the Lusternik-Schnirelmann category of $X$.

On the other hand, we denote by  ${\mathbb L}(S)$  the free graded Lie algebra on a set $S= S_0$. If $I\subset \mathbb L(S)$ is an ideal generated by elements homogeneous with respect to the length of iterated brackets in $S$, then
$$\mathbb L(S)/I = \oplus_n \, \mathbb L(S)(n)$$
identifies $\mathbb L(S)/I$ as a Lie algebra with a weight decomposition, where $\mathbb L(S)(n)$ is the image of the Lie brackets of length $n$ in elements of $S$. 

We can now state our two main theorems. 
 
 Recall first that the fundamental group $G$ of a closed orientable Riemann surface $X_g$ of genus $g>0$ is the free group on generators $a_i$, $b_i$, $1\leq i\leq g$, divided by the single relation $\prod_{i=1}^g a_ib_ia_i^{-1}b_i^{-1}$. Associated with $X_g$ is the weighted Lie algebra, 
 $$L := \mathbb L(a_i, b_i)_{1\leq i\leq g}/ \left( \, \sum_{i=1}^g[a_i, b_i]\,\right)\,.$$

 In \cite{La}, Labute proves that the weighted Lie algebras $L$ and $L_G$ are isomorphic. We recover this result and we prove:

\vspace{3mm}\noindent {\bf Theorem 1.}  {\sl  Let $  X_g$ be a     Riemann surface   of genus $g>0$ with fundamental group $G$, associated Lie algebra $L$ and minimal Sullivan model $(\land V,d)$.  Then
\begin{enumerate}
\item[(i)] The weighted Lie algebras $L_G$ and $L$ are isomorphic.   
\item[(ii)]  $$\mbox{depth}\, (\land V,d) = \mbox{depth}\, L= \mbox{depth}\, UL= \mbox{depth}\, \mathbb Q[G] = 2\,.$$
\end{enumerate}}

Next recall that a right-angled Artin group is a group $A$ that admits a presentation of the form
$$A = <\, \{ x_i\}_{1\leq i\leq n};  x_ix_j  =  x_jx_i \,, \hspace{3mm}\mbox{for some possibly void set ${\mathcal S}$ of pairs } (i,j)>\,.$$
Right-angled Artin groups include finitely generated free abelian groups   and   finitely generated free groups.  

Associated to $A$ is a flag polyhedron $P$ whose $r$-simplices are the subsets $\sigma = <x_{i_0}, \cdots ,x_{i_r}>$ of generators in which the $x_{i_j}$ all commute. The centralizer $Z(\sigma)\supset \sigma$ is the subset of generators which commute with each of the $x_{i_j}$. We say $\sigma$ is \emph{disconnecting} if the other generators in $Z(\sigma)$ can be divided into two sets $\{y_i\}$ and $\{z_j\}$ such that for all $i,j$, $y_iz_j\neq z_jy_i$.  

Also associated to $A$ is the weighted  Lie algebra $L$ defined by
$$L = \mathbb L(x_i)_{ 1\leq i\leq n}/I\,,$$
where $I$ is the ideal generated by the brackets $[x_i,x_j]$ with $(i,j)\in \mathcal S$. In \cite{PS}, Papadima and Suciu show that $L$ and $L_A$ are isomorphic. We recover this fact, and  building on a theorem of Jensen and Meier \cite{JM} we have

 \vspace{3mm}\noindent {\bf Theorem 2.} {\sl Let $A$ be a right-angled Artin group with associated Lie algebra $L$, and let $(\land V,d)$ be the minimal Sullivan model of a classifying space $K$ for $A$. Then
 \begin{enumerate}
 \item[(i)] The weighted Lie algebras $L_A$ and $L$ are isomorphic.
 \item[(ii)] 
 $$ \mbox{depth}\,  (\land V,d)= \mbox{depth}\, L= \mbox{depth}\, UL=   \mbox{depth}\, \mathbb Q[A]\,.$$
\item[(iii)] If $A$ is abelian then depth$\, \mathbb Q[A]= n$. Otherwise depth$\, \mathbb Q[A]=$ least $r+2$ for which there is a possibly empty disconnecting $r$-simplex $\sigma$. 
 \end{enumerate}}

\vspace{3mm} Theorem 2  permits the computation of depth $\mathbb Q[A]$ in non trivial examples.
For instance, suppose the generators of $A$ are divided into three groups, $\{x_i\}, \{y_j\}$ and $\{z_k\}$, in which the $x_i$ are central and $y_jz_k\neq z_ky_j$ for all $j$ and $k$. If there are $r$ generators $x_i$ then 
$$\mbox{depth}\, \mathbb Q[A] = r+1\,.$$

On the other hand, suppose $A$ has generators $x_i, y_j$ in which, in addition to relations among the $x_i$ and among the $y_j$ there is a single relation $x_1y_1=y_1x_1$. Then
$$\mbox{depth}\, \mathbb Q[A]\leq 3\,.$$
with the exact value depending on the other relations.

\vspace{3mm}  
  The paper is organized as follows. In section 2 we  review Sullivan models $(\land V,d)$ and define depth$\, (\land V,d)$. We recall the definition of the homotopy Lie algebra $L_V$ of $(\land V,d)$ and the weighted Lie algebra 
$$E_V = \bigoplus_n E_V(n)\,, \hspace{1cm} E_V(n) := (L_V^n)_0 \, /\, (L_V^{n+1})_0\,.$$
If a  path connected space $X$, with fundamental group $G$ and minimal Sullivan model, $(\land V,d)$, satisfies   dim$\,H^1(X) <\infty$, (\cite[Theorem 7.7]{FHTII}) then there is a natural isomorphism of weighted Lie algebras
\begin{eqnarray}
\label{UUn}
L_G \cong E_V\,.\end{eqnarray}

In section 3 we consider Sullivan algebras $(\land V,d)$ satisfying the conditions
$$(\land V,d) \mbox{ is formal},\hspace{3mm}  V= V^1, \mbox{ and dim}\, H(\land V,d)<\infty\,.$$
For these Sullivan algebras we show that
$$\mbox{depth}\, (\land V,d)= \mbox{depth}\, E_V\,.$$
Then in section 4 we show that for these Sullivan algebras   
$$\mbox{depth}\, E_V= \mbox{ depth}\, UE_V\,.$$
 
Section 5 is devoted to the example of orientable closed surfaces. Section 6  sets up the machinery for right-angled Artin groups, $A$, and section 7 contains the proof of Theorem 2. This relies  on a combinatorial description of  depth$\, \mathbb Q[A]$ established in \cite{JM}.

\section{Sullivan models} 

  Denote by $\land V$   the free graded commutative algebra on a graded vector space $V$, and by $\land^kV$   the linear span of the monomials of length $k$.

\vspace{3mm}\noindent {\bf Definition.} A \emph{Sullivan algebra} is a commutative differential graded algebra (cdga for short) of the form $(\land V,d)$,  where   $V= V^{\geq 1}$  is equipped with a filtration $V = \cup_{p\geq 0} V_p$, with $d(V_0)=0$, 
$V_{p-1}\subset V_p$, and $$ d(V_p) \subset \land (V_{p-1})\,,\hspace{3mm}\mbox{ $p\geq 1$}\,.$$ The Sullivan algebra $(\land V,d)$ is called \emph{minimal} if $$d(V) \subset \land^{\geq 2}V\,.$$

To each path connected space $X$, Sullivan   associates  the cdga, $A_{PL}(X)$,  of rational polynomial forms on the simplicial set Sing$\,X$, of singular simplices on $X$. This is a rational analogue of the cdga of de Rham forms on a manifold. There is then a (unique up to isomorphism) minimal Sullivan algebra $(\land V,d)$ together with a cdga quasi-isomorphism $$\varphi : \xymatrix{(\land V,d)\ar[rr] && A_{PL}(X)}\,.$$
This is the \emph{minimal Sullivan model of $X$}.

For any minimal Sullivan algebra, $(\land V,d)$, denote by $d_1$ the quadratic part of the differential $d$. It is   characterized by the properties
$$d_1 : V\to \land^2V \hspace{5mm}\mbox{and }  d-d_1 : V\to \land^{>2}V\,.$$
  The \emph{homotopy Lie algebra} of $(\land V,d)$ is the graded Lie algebra  $L_V = \{(L_V)_p\}_{ p\geq 0}$, where
  $$(L_V)_p :=\mbox{Hom}(V^{p+1}, \mathbb Q)\,.$$ 
  The associated pairing $V\otimes sL_V \to \mathbb Q, (v\otimes f)\mapsto <v,f>:=(-1)^{deg\, v} f(v)$ extends to the pairing between 
  $\land^2V$ and $\land^2sL_V$ given by
  $$<v_1\land v_2, f\land g> := <v_1, g>\cdot <v_2, f> + (-1)^{deg\, f\cdot deg\, g} <v_1, f>\cdot <v_2,g>\,.$$
 The Lie bracket in $L_V$ is then defined by 
$$<v, s[x,y]> := (-1)^{deg\, y\,+1} <d_1v, sx,sy>\,.$$

Next, denote by $(L_V^n)_0$ the ideal generated by Lie brackets of length $n$ in $(L_V)_0$. Then a weighted Lie algebra, $E_V$, is given by
$$E_V= \bigoplus_{n\geq 1} E_V(n)\,, \hspace{1cm} E_V(n) = (L_V^n)_0\,/\, (L_V^{n+1})_0\,.$$
Here, in analogy with the rational Malcev Lie algebra, the Lie bracket in $E_V$ is induced by the commutator
$$[\,\,,\,\,]: (L_V^p)_0 \times (L_V^q)_0\to (L_V^{p+q})_0\,.$$

\vspace{2mm} Another important construction in Sullivan theory is the acyclic closure, $(\land V\otimes \land U,d)$, of a minimal Sullivan algebra $(\land V,d)$. It is a cdga quasi-isomorphic to $(\mathbb Q, 0)$, which  extends $(\land V,d)$, and    satisfies $d(U) \subset \land^+V\otimes \land U$.   It satisfies the same filtration condition as $(\land V,d)$; however, $U^0$ will be non-zero if $V^1\neq 0$.   

Denote by $d_1(1\otimes \Phi)$ the component of $d(1\otimes \Phi)$ in $V\otimes \land U$. Then the \emph{holonomy representation} $\theta$ of $L_V$ in $\land U$ is defined by 
$$\theta (x)\Phi := -<d_1(1\otimes \Phi),x>\,.$$

\vspace{3mm} Let $(\land V\otimes \land U, d)$ be the acyclic closure of $(\land V,d)$. Then $(\land V\otimes \land U, d_1)$ is the acyclic closure of $(\land V, d_1)$. Denote by Hom$^p_{\land V}(\land V\otimes \land U, \land V)\subset \mbox{Hom}_{\land V}(\land V\otimes \land U, \land V)$ the subspace of maps $\land^rV\otimes \land U\to \land^{p+r}V$, $r\geq 0$.   Then the differential $d_1\circ f- (-1)^{deg\, f} f\circ d_1$ yields a complex
$$\dots \to \mbox{Hom}^p_{\land V}(\land V\otimes \land U, \land V) \to \mbox{Hom}_{\land V}^{p+1}(\land V\otimes \land U, \land V)\to \dots$$
with homology   Ext$^p_{(\land V,d_1)}(\mathbb Q, (\land V,d_1))$.  

The \emph{depth} of a minimal Sullivan algebra (\cite[Section 10.2]{FHTII}) is then defined by
 $$\mbox{depth$\, (\land V,d):= $ least $p$ (or $\infty$) such that Ext}^p_{(\land V,d_1)}(\mathbb Q, (\land V, d_1))\neq 0\,.$$

\section{Formal Sullivan algebras}

A Sullivan algebra $(\land V,d)$ is \emph{formal} if there is a quasi-isomorphism $(\land V,d) \to H(\land V,d)$, in which case this may be chosen to induce the identity in $H(\land V,d)$. As will be shown in sections 6 and 7, the minimal Sullivan models, $(\land V,d)$, of closed orientable Riemann surfaces, and of the classifying spaces of right-angled Artin groups, have the following three properties:
\begin{eqnarray}
\label{duo}
(\land V,d) \mbox{ is formal, $V = V^1$, and dim}\, H(\land V,d)<\infty\,.
\end{eqnarray}

\vspace{3mm}\noindent {\bf Proposition 1.}  (\cite[\S 15.3]{FHTII}).   {\sl Suppose $(\land V,d)$ is a minimal Sullivan algebra satisfying (\ref{duo}). Then 
 \begin{enumerate}
 \item[(i)] $V =\oplus_{n\geq 0} V(n)$ with $d : V(0)\to 0$ and  $d: V(n)\to +_{p+q=n} V(p)\land V(q)$.
 \item[(ii)]   dim$\, V(n)<\infty$, $n\geq 0$.
 \item[(iii)] The inclusion of $V(0)$ in $V$ induces an isomorphism $V(0)\stackrel{\cong}{\rightarrow} H^1(\land V,d)$.
 \item[(iv)] $(\land V,d)$ is the direct sum of the finite dimensional subcomplexes:
 $$C_V(r): = \bigoplus_{p+n=r} ((\land^p V)(n),d)\,,$$
 where  $(\land^pV)(n) = +_{n_1+\dots + n_p= n}\, V({n_1})\land \dots \land V({n_p})$.
 \end{enumerate}
}

\vspace{3mm}
From the definition of the Lie bracket in $L_V$,
 $$E_V  := \bigoplus_n \mbox{Hom}(sV(n), \mathbb Q) \subset L_V$$
is a sub Lie algebra. Moreover,  \cite[Theorem 2.1]{FHTII} gives
$$L_V^r = \mbox{Hom}(\bigoplus_{n\geq r-1} sV(n), \mathbb Q)\,.$$
This identifies $E_V$ as the associated bigraded Lie algebra for $L_V$, with
$$E_V(n):= \mbox{Hom}(sV(n-1), \mathbb Q) = L_V^n/L_V^{n+1}\,.$$ Note that   $E_V(n)$ is the linear span of the Lie brackets of length $n$ in $E_V(1)$.

\vspace{3mm}\noindent {\bf Proposition 2.} {\sl Suppose $(\land V,d)$ is a minimal Sullivan algebra satisfying (\ref{duo}).  Then
\begin{enumerate}
\item[(i)] depth$\, (\land V,d) =$ depth$\, E_V$.
\item[(ii)] $H(C_*(E_V))$ is finite dimensional, where $C_*(E_V) = \land (sE_V, \partial_1)$ is the classical Cartan-Chevalley-Eilenberg construction.
\end{enumerate}
}

\vspace{3mm}\noindent {\sl proof.} (i) It is immediate from the construction that 
$$(\land V,d) = \varinjlim_n C^*(L_V/L_V^n)= \varinjlim_n C^*(E_V/E_V^n)\,.$$
Thus by \cite[Theorem 15.1]{FHTII},
$$\mbox{depth}\, (\land V,d)= \mbox{depth}\, E_V\,.$$

\vspace{2mm} (ii) Since  dim$\,H^*(\land V,d)<\infty$, there is some $r_0$ such that   $H^*(\, C_V(r)\,)= 0$, $r>r_0$. On the other hand it is immediate from the definition  that $(\land sE, \partial_1)$ is the direct sum of   subcomplexes dual to the finite dimensional complexes $(C_V(r),d)$, and so dim $\,H(\land sE_V, \partial_1)<\infty$. 
  \hfill$\square$
  
 \vspace{3mm}Finally, let $S$ be a copy of $(E_V)_1$, and extend the identification to a surjection
 $$\pi: \mathbb L(S)\to  E_V\,.$$
 Denote by $I$ the ideal generated by  $\pi(2) : [S,S]\to (E_V)(2)$. Then \cite[Theorem 15.4]{FHTII}, together with the discussion above,  yields
 
 \vspace{3mm}\noindent {\bf Proposition 3.} (see also \cite{Chen}). {\sl The surjection $\pi$ factors to yield   an isomorphism of weighted Lie algebras,
 $$\mathbb L(S)/I \stackrel{\cong}{\longrightarrow} E_V\,.$$  }

\section{Depth of a weighted Lie algebra}

 In this section we consider  graded Lie algebras, $L$,  of the form $\mathbb L(S)/I$, where $S = S_0$ and $I$ is generated by elements that are homogeneous with respect to the length of iterated Lie brackets in $S$.  
 
 Denote by $\mathbb L(r)$ the linear span of the iterated Lie brackets of length $r$ in $S$, and by $L(r)$ the image in $L$ of $\mathbb L(r)$. By hypothesis, $L$ is a weighted Lie algebra:
$$L = \bigoplus_{r\geq 1} L(r)$$
and $[L(r), L(s)] \subset L(r+s)]$. Moreover $L^{r+1} = \oplus_{s\geq r+1} L(s)$ and so 
$$ \bigoplus_{s\leq r} L(s)\stackrel{\cong}{\longrightarrow} L/L^{r+1} \,.$$

\vspace{3mm}\noindent {\bf Remark.} Suppose $(\land V,d)$ is a Sullivan algebra satisfying (\ref{duo}). Then Proposition 3 asserts that the weighted Lie algebra $L= E_V$ associated with $L_V$ has the form above, and that   the subspace $E_V(r)$ defined in the previous section coincides with $L(r)$ as defined above.

\vspace{3mm}\noindent {\bf Proposition 4.} With the  notation and hypotheses above, suppose that dim$\, H(C_*(L))<\infty$. Then
$$\mbox{depth}\, UL= \mbox{depth}\, L=\mbox{depth}\, \widehat{UL}\,.$$

\vspace{3mm}\noindent {\sl proof.} 
The direct decomposition of $L$ extends to direct decompositions,
$$UL = \bigoplus_{r\geq 0} (UL)(r) \hspace{5mm}\mbox{and}\hspace{2mm} \land sL = \bigoplus_{p,r\geq 0} (\land^psL)(r)\,,$$
characterized by the properties $(UL)(0)= \mathbb Q=(\land sL)(0)$,
$$(UL)(r)\cdot (UL)(s) \subset (UL)(r+s), \hspace{5mm}\mbox{and}\hspace{2mm} (\land^p sL)(r)\land (\land^qsL)(s) \subset (\land^{p+q}sL)(r+s)\,.$$
Moreover in the Cartan-Chevalley-Eilenberg construction, $(C,\partial) = (\land sL\otimes UL, \partial)$, $$ \partial = \partial_2\otimes id + \partial_1\,,$$
$$\partial_2 : (\land^psL)(r) \to (\land^{p-1}sL)(r),\hspace{4mm}\mbox{and}\hspace{2mm} 
\partial_1 : (\land^psL)(r)\otimes 1 \to \bigoplus_{q<r} (\land^{p-1}sL)(q) \otimes UL(r-q)\,.$$

By hypothesis, $H(C_*(L)) = H(\land sL, \partial_2)$  is finite dimensional. But $H(\land sL) = \oplus_r H((\land sL)(r))$. It follows that for some $k$,
$$H((\land sL)(r), \partial_2) = 0\,, \hspace{1cm} r>k\,.$$
Set $R = \oplus_{r\leq k} (\land sL)(r) $. Then $H(R) \stackrel{\cong}{\rightarrow} H(\land sL, \partial_2)$. Moreover $R\otimes UL$   is preserved by both $\partial_2\otimes id$ and $\partial_1$. 

Recall now that $(\land sL\otimes UL,\partial)$ is a free resolution of $\mathbb Q$ by $UL$-modules. An easy spectral sequence argument then shows that the inclusion
$$(R\otimes UL), \partial)\to (\land sL\otimes UL, \partial)$$
is a quasi-isomorphism. Thus writing $R_p = R\cap \land^psL$ we have that
$$\cdots \stackrel{\partial}{\longrightarrow} R_{p-1}\otimes UL \stackrel{\partial}{\longrightarrow} R_p\otimes UL \stackrel{\partial}{\longrightarrow} \cdots$$
is a $UL$-free resolution of $\mathbb Q$. Moreover, since each $(\land sL)(r)$ is finite dimensional $R$ itself is finite dimensional.

Next observe that since $UL$ is generated by $L(1)$ it follows that the augmentation ideal $I\subset UL$ satisfies
$$I^{n+1} = \oplus_{r>n} (UL)(r)\,.$$
This identifies the inclusion $\lambda : UL \to \widehat{UL}$ as the inclusion $\bigoplus_{r\geq 0} UL(r) \to \prod_{r\geq 0} UL(r)$. In particular, since dim$\, R<\infty$, this identifies
$$\mbox{Hom}_{UL}(R\otimes UL, UL) \to \mbox{Hom}_{UL}(R\otimes UL, \widehat{UL})$$
 as the inclusion
 $$\bigoplus_n \left( \bigoplus_{q+i=n} \mbox{Hom} (R(q), (UL)(i))\right) \to \prod_n \left(\bigoplus_{q+i=n} \mbox{Hom}(R(q), (UL)(i))\right)\,.$$
 Since each $\oplus_{q+i=n} \mbox{Hom}(R(q), (UL)(i))$ is a subcomplex, 
it follows that
 $$\renewcommand{\arraystretch}{1.6}
 \begin{array}{ll}
 \mbox{depth}\, UL &=  \mbox{least $p$ (or $\infty$) such that Ext}^p_{UL}(\mathbb Q,   {UL})\neq 0\,,\\
 &= \mbox{least $p$ (or $\infty$) such that Ext}^p_{UL}(\mathbb Q,  \widehat{UL})\neq 0\,,\\
 & = \mbox{depth}\, L\,.
 \end{array}
 \renewcommand{\arraystretch}{1}
 $$
 
 On the other hand, $(R\otimes UL, \partial)$ is the direct sum of the subcomplexes $$C(n):= \bigoplus_{q+i=n} R(q)\otimes (UL)(i)\,.$$ It follows that unless $n= 0$, 
 $H(C(n)) = 0\,.$ 
 Since $R\otimes \widehat{UL}$ is given by
 $$R\otimes\widehat{UL} = \prod_n \left[ \bigoplus_{q+i=n} R(q)\otimes (UL)(i)\right]$$
 it follows that $(R\otimes \widehat{UL}, \partial)$ is a $\widehat{UL}$-free resolution of $\mathbb Q$. Since
 $$\mbox{Hom}_{\widehat{UL}}(R\otimes \widehat{UL}, \widehat{UL}) = \mbox{Hom}_{UL}(R\otimes UL, \widehat{UL}),$$
this yields  
 $$\mbox{depth}\, \widehat{UL} = \mbox{depth}\, L\,.$$
 
 \hfill$\square$

 \vspace{5mm}\noindent {\bf Proposition 5.} {\sl Suppose $F = \oplus_r F(r)$ is a sub weighted Lie algebra of $L$ for which dim$\, H(C_*(F))<\infty$. Then
$$\mbox{depth}\, UF= \mbox{least $p$ (or $\infty$) such that Ext}^p_{UF}(\mathbb Q, UL)\neq 0\,.$$}

\vspace{3mm}\noindent {\sl proof.} As with $UL$ we have a free resolution of $UF$ of the form
$$R_F\otimes UF\stackrel{\simeq}{\longrightarrow} \mathbb Q$$
in which dim$\,R_F<\infty$. Then
$$\mbox{Hom}_{UF}(R_F\otimes UF, UL)= \mbox{Hom}(R_F, UL)\,.$$
Since $UL$ is a free $UF$-module we may write $UL = UF\otimes W$ and 
$$\mbox{Hom}(R_F, UL) = \mbox{Hom}(R_F, UF)\otimes W.$$ The Proposition follows. \hfill$\square$

\section{Depth of a Riemann surface}

Let $X_g$ be an orientable Riemann surface of genus $g\geq 1$ with fundamental group $G$.  The space $X_g$ is   formal (\cite{DGMS}), dim$\, H(X_g)<\infty$,   and the minimal Sullivan model
$(\land V,d)$ satisfies $V= V^1$ (\cite[Theorem 8.1]{FHTII}). Thus $(\land V,d)$ satisfies (\ref{duo}). In particular, by Proposition 1,  $$V = \oplus_{n\geq 0} V(n)\hspace{3mm}\mbox{with } d(V(n)) \subset +_{p+q=n} V(p)\land V(q)\,.$$

 The components $V(0)$ and $V(1)$ admit the following explicit descriptions, as described in \cite[\S 8.5]{FHTII}.
A basis for $V(0)$ is given by the elements  $u_i, v_i$ with $1\leq i\leq g$.
A basis of 
$V(1)$ is given by   elements  $u_{ij}, v_{ij} $ with $1\leq i<j\leq g$, together with elements $w_{k\ell}$, and $1\leq k,\ell\leq g$ with $k+\ell >2$.
The differential $d : V(1)\to \land^2V(0)$ is given by
$$du_{ij} = u_iu_j\,, \hspace{5mm}
 dv_{ij} = v_iv_j\,, \hspace{5mm} 
 dw_{k\ell} = u_kv_\ell \hspace{3mm} \mbox{for }k\neq \ell\,, $$
 $$\mbox{and }
 dw_{kk} = u_kv_k- u_1v_1\,.$$

Next, let $a_i, b_j$ ($1\leq i,j\leq g$) be the basis of $E_V(1)$ dual to the elements $su_i$, $sv_j$. Then the Lie brackets $[a_i, a_j]$ and $[b_i, b_j]$, $i<j$ are respectively dual to $su_{ij}$ and $sv_{ij}$. Moreover, for $i\neq j$, $[a_i, b_j]$ is dual to $sw_{ij}$. Finally, for $k>1$, $[a_k, b_k]-[a_1, b_1]$ is dual to $sw_{kk}$, while
$$<\sum\, [a_i, b_i], sE_V(2)> = 0\,.$$
It follows that ker$\, \mathbb L(a_i, b_i)(2)\to E_V(2)$ is just $\mathbb Q\cdot \sum[a_i,b_i]$. Thus Proposition 3 shows that the surjection $\mathbb L(a_i, b_i) \to E_V$ induces an isomorphism
$$L:= \mathbb L(a_i, b_i)/\sum[a_i, b_i] \stackrel{\cong}{\longrightarrow} E_V\,.$$


\vspace{3mm}\noindent {\bf Theorem 1.} Let $X_g$ be an orientable closed Riemann surface of genus $g>0$, with fundamental group $G$, and minimal Sullivan model $(\land V,d)$. Then with the   notation above,
\begin{enumerate}
\item[(i)] The weighted Lie algebras $L_G$ and $L$ are isomorphic.
\item[(ii)] $$\mbox{depth}\, (\land V,d) = \mbox{depth}\, L= \mbox{depth}\, UL= \mbox{depth}\, \mathbb Q[G]= 2\,.$$
\end{enumerate}

\vspace{3mm}\noindent {\sl proof.} (i). By (\ref{UUn}) in the Introduction,   $L_G\cong E_V$, and, as observed above, $E_V\cong L$.

(ii)  The group $G$ is a Poincar\'e duality group and it follows from \cite[Chap VIII, Proposition 8.2]{Br} that depth$\, \mathbb Q[G]= 2$.   On the other hand, by Proposition 2 and 3, depth$\, (\land V,d)=$ depth$\, E_V$ and depth$\, L=$ depth$\, UL$. Thus it remains to show that depth$\, (\land V,d)= 2$.

Write $H = H(\land V,d)$. By definition there is a quasi-isomorphism $\varphi : (\land V,d) \to (H,0)$. Since $V = V^1$, $d= d_1$ and $\varphi (\land^pV)\subset H^p$. Therefore $\varphi$ induces isomorphisms
$$\mbox{Ext}^p_{\land V}(\mathbb Q, \land V)\stackrel{\cong}{\longrightarrow} \mbox{Ext}^p_{\land V}(\mathbb Q, H)\stackrel{\cong}{\longleftarrow} \mbox{Ext}^p_H(\mathbb Q, H)\,.$$

On the other hand, since $H$ satisfies Poincar\'e duality, $\mbox{Hom}(H, \mathbb Q)$ is a free $H$-module of rank one via the cap product. It follows that $$\mbox{Ext}_H(\mathbb Q, H) = \mbox{Hom}\left(\, \mbox{Tor}^H(\mathbb Q, \mbox{Hom}(H,\mathbb Q), \mathbb Q\, \right)$$
has dimension one. Finally let $\omega \in H^2$ denote the fundamental class. Then the map $f\in  \mbox{Hom}_H(\mathbb Q, H)$ defined by $f(1) = \omega$ is a non trivial cycle.  Since Ext$_H(\mathbb Q,H)$ has dimension 1, Ext$_H(\mathbb Q, H) = \mathbb Q\cdot f$. Therefore depth$\, (\land V,d)=2$. \hfill$\square$

 \vspace{4mm}Now consider the completion $\widehat{L}$ of the weighted Lie algebra $L$,
 $$\widehat{L}= \varprojlim_n L/L^n\,.$$
 This Lie algebra $\widehat{L}$ is the rational homotopy Lie algebra of the minimal model $(\land V,d)$.
 
 \vspace{3mm}\noindent {\bf Proposition 6.} {\sl Let $L$ be the weighted Lie algebra associated to $X_g$. Then
 $$\mbox{depth}\, \widehat{L}=\mbox{depth}\, L = 2\,.$$}
 
 \vspace{3mm}\noindent {\sl proof.} Since $\widehat{L}$ has no zero divisor, depth$\, \widehat{L} \geq 1$.  On the other hand by 
 \cite[Theorem 2]{Sdep}, depth$\, \widehat{L}\leq $ cat$(\land V,d)= 2$. It remains to show that 
 Ext$^1_{U\widehat{L}}(l\!k, \widehat{U\widehat{L}})= 0$. 

 First observe from \cite[Lemma 2.12]{FHTII} that $\widehat{UL}= \widehat{U\widehat{L}}$. 

Suppose $\omega \in \mbox{Ext}^1_{U\widehat{L}}(l\!k, \widehat{UL})$. By hypothesis, $\omega$ may be represented by a $U\widehat{L}$-linear map
 $$f :  s\widehat{L}\otimes U\widehat{L} \to \widehat{UL}$$
 such that $f\circ d= 0$ and $f(sL\otimes 1) = 0$. Since depth$\, L=2$, we may assume that $f$ vanishes on $sL\otimes U\widehat{L}$. 

Assume by induction   that $f : s\widehat{L}\to I^n$, some $n\geq 0$, where $I$ is the augmentation ideal in $U\widehat{L}$. Then according to \cite[Theorem 2.1]{FHTII}, $\widehat{L} = L+ [L, \widehat{L}]$. If $x\in s\widehat{L}$ then $x= x_0+ \sum [x_i, y_i]$ and, because $d\circ f= 0$ and $f(sL)= 0$,
$$f(sx) = f(\sum s[x_i, y_i]) = \sum\pm f(sx_i)y_i \pm f(sy_i)x_i \in I^{n+1}\,.$$
It follows that $f(s\widehat{L})\in \cap_n I^n = 0$. \hfill$\square$

 \section{Right-angled Artin groups}
 
 Right-angled Artin groups are particular Artin groups,  introduced by A. Baudisch in the 70's (\cite{Bau}), and further developed in the 80's by C. Droms under the name of graph groups. They have been the subject of considerable research since that time, with a   survey  by R. Charney in 2007 \cite{Char}. Recall the definition from the Introduction:
 
\vspace{3mm}\noindent {\bf Definition.} 
  A right-angled Artin group $A$ is a group with presentation of the form
 $$A = <x_1, \ldots, x_n \,\vert\,   x_ix_j  =   x_jx_i    \,\, \mbox{for some possibly void subset $\mathcal S$ of pairs $(i,j)$}>\,.$$

{\sl For the rest of this section we fix a right-angled Artin group, $A$, together with a presentation of the form above.}

\vspace{3mm} With this presentation of $A$ are associated:
\begin{enumerate}
\item[$\bullet$] A \emph{graph} $\Gamma$. The vertices are labelled by the generators $x_i$ and   a pair $x_i, x_j$ is connected by an edge if  and only if $x_ix_j= x_jx_i$.
\item[$\bullet$] The \emph{associated flag polyhedron}, $P$,  the polyhedron with vertices the $x_i$ and $(r-1)$ simplices $<x_{i_1}, \cdots , x_{i_r}>$, $i_1<\cdots <i_r$, corresponding to the complete subgraphs of $\Gamma$.  The empty simplex is denoted by $<\emptyset>$. 
\item[$\bullet$] A \emph{commutative graded algebra} $C = \land (t_1, \ldots , t_n)/ I$, where each deg$\, t_i= 1$ and $I$ is generated by the monomials $t_i\land t_j$ for which $x_ix_j\neq x_jx_i$. This  is   the \emph{flag algebra} for $\Gamma$.  A monomial $t_{i_1}\land \cdots \land  t_{i_r}$ in $C$ is nonzero if and only if $x_{i_1}, \cdots, x_{i_r}$ are the vertices of a   complete sub graph of $\Gamma$, and those monomials form a basis for $C$.

\item[$\bullet$] The \emph{dual coalgebra} $B = \mbox{Hom}(C,\mathbb Q)$. We denote by $sx_1, \cdots , sx_n$ the basis of $B_1$ dual to the basis $t_1, \cdots , t_n$ of $C^1$: $<t_i, sx_j>= \delta_{ij}$. Thus the basis of $B_r$ dual to the basis of $C^r$ consists of the elements $ <sx_{i_1}, \cdots , sx_{i_r}>$, $i_1<\cdots <i_r$, where $x_{i_1}, \cdots x_{i_r}$ are the vertices of some complete subgraph of $\Gamma$.  This basis of $B_r$ is identified with the set of $(r-1)$-simplices of $P$ via $<sx_{i_1}, \cdots , sx_{i_r}> \longleftrightarrow <x_{i_1}, \cdots , x_{i_r}>$ and we will abuse (and simplify) notation by using $<x_{i_1}, \cdots , x_{i_r}>$ to denote both.  Finally, we denote by $e=<\emptyset>$ the element of $B_0$ dual to $1\in C_0$.

\item[$\bullet$] A \emph{CW complex} $K$, called the \emph{Salvetti complex},   the union of tori
$$K = \cup_{<x_{i_1}, \cdots , x_{i_r}>\in P}\,\,\, S^1_{i_1}\times S^1_{i_2}\times \dots \times S^1_{i_r}\,.$$
\item[$\bullet$] The \emph{minimal Sullivan model}, $(\land V,d)$ for $K$. 
\item[$\bullet$]  The  \emph{weighted Lie algebra} $L$ of Theorem 2: $$L= \mathbb L(x_i)_{1\leq i\leq n}/I\,,$$
where $I$ is   the ideal generated by the $[x_i, x_j]$ for the couples $(i,j)$ corresponding to edges in $P$.
\end{enumerate}

\vspace{3mm} Immediately from the definition is the well known

\vspace{3mm}\noindent {\bf Remark.}\begin{enumerate}
\item[(i)] Every finite graph is the graph of a unique right-angled Artin group.
\item[(ii)] A finite polyhedron $P'$ is the polyhedron of a right-angled Artin group if and only if whenever the faces of a simplex of dimension $\geq 2$  are in $P'$ then the simplex is also in $P'$.  
\end{enumerate}

\vspace{3mm} The relations between   the algebraic invariants described above are reflected in  the Fr\"oberg resolution:

\vspace{3mm}\noindent {\bf Proposition 7.} (\cite{Fr}) {\sl  A resolution of $\mathbb Q$ by free right $UL$-modules is defined by the complex,  
  $$\cdots \to   B_r \otimes UL   \stackrel{d}{\longrightarrow}   B_{r-1} \otimes UL  \to \cdots \to B_1\otimes UL\to   B_0\otimes UL  \stackrel{\varepsilon}{\to} \mathbb Q\,,$$
where  $\varepsilon: UL\to \mathbb Q$ is the augmentation, $d  ( <x>) = e\otimes x$, and for $r\geq 2$, 
  $$d( <x_{i_1}, \cdots , x_{i_r}>) = \sum_{j=1}^{r} (-1)^{j-1}    <x_{i_1}, \cdots ,\widehat{x_{i_j}} ,  \cdots , x_{i_r}>\otimes x_{i_j} \,.$$
}

\vspace{3mm}  The homotopy type of the Salvetti complex has been described by Charney and Davis:

\vspace{3mm}\noindent {\bf Theorem 3.}  (\cite{CD}) The Salvetti complex $K$ is a $K(A,1)$-space.

\vspace{3mm} It follows from the work of  Kapovich and Millson   (\cite{Mis}), and Papadima and Suciu (\cite{PS}) that the Salvetti complex is a formal space and that its minimal Sullivan model satisfies    $V= V^1$. We give a direct proof of those facts.

\vspace{3mm}\noindent {\bf Theorem 4.}  {\sl  With the definitions and notation above:
\begin{enumerate}
\item[(i)]  The Salvetti complex $K $ is a formal space;
\item[(ii)] $H^*(K;\mathbb Q)$ is isomorphic as a graded algebra to the algebra $C$;
\item[(iii)] The minimal Sullivan model $(\land V,d)$ of $K$ satisfies $V = V^1$;
\end{enumerate}
In particular, $(\land V,d)$ satisfies (2).}

\vspace{3mm}\noindent {\sl proof.} (i).     The Salvetti complex $K$ is a  union of tori $T_i$. We can thus write $K = \cup_{p\leq r } X_p$ where $\{*\}= X_0\subset X_1\subset X_2\subset \cdots \subset X_r= K$, with the following properties: each $X_p$ is a union of   tori,  and for $p\geq 0$,  $H_*(X_p)\to H_*(X_{p+1})$ is injective with cokernel of dimension one. By induction we can suppose that $X_p$ is formal and that a model for the injection of each sub complex  is given by the induced morphism in cohomology. 

Given a $k$-torus $T = S^1_1\times \cdots \times S^1_k$, denote by $\partial T$ the union of the codimension one subtori: $\partial T$ is the $(k-1)$-skeleton in the standard CW structure. We may suppose that $X_{p+1}$ is obtained by adding to $X_p$ a torus $T$  such that $\partial T$ is a subcomplex of $X_p$. Since models for the injection $\partial T\to T$ and $\partial T\hookrightarrow X_p$ are given by the cohomology maps, it follows from (\cite[Proposition 13.6]{FHT}) that $X_{p+1}$ is formal and its cohomology is given by the pullback $H^*(T)\times_{H^*(\partial T)} H^*(X_p)$. 

\vspace{2mm} (ii) The Salvetti complex is a subcomplex of $S^1_1\times \cdots \times S^1_n$. The canonical isomorphism $\land (t_1, \cdots , t_n)\to H^*(S^1_1\times \cdots \times S^1_n;\mathbb Q)$ factors to give a morphism $C = \land (t_1, \cdots , t_n)/I\to H^*(K;\mathbb Q)$. The same limit argument as above shows that this is an isomorphism.

\vspace{2mm} (iii). Observe that, as in Section 5,
$$L=  \bigoplus_{r\geq 1} L(r)\hspace{5mm}\mbox{and } UL= \bigoplus_{r\geq 0} UL(r)\,,$$
where $L(r)$ is the linear span of the iterated Lie brackets of length $r$ in the generators $x_i$ and $UL$ is the linear span of the monomials of length $r$ in the $x_i$. 

As the dual of $C$, $B$ is a coalgebra with comultiplication  given by
$$\renewcommand{\arraystretch}{1.6}
\begin{array}{ll} \Delta ( <x_{i_1}, \cdots ,x_{i_r}>) &= <x_{i_1}, \cdots, x_{i_r}> \otimes 1 \\&+ \sum_{(p, \lambda)\in T}(-1)^\lambda <x_{i_{\lambda (1)}}, \cdots , x_{i_{\lambda (p)}}> \otimes <x_{i_{\lambda(p+1)}}, \cdots , x_{i_{\lambda r)}}> \\ &
+ 1\otimes <x_{i_1}, \cdots , x_{i_r}>\,.
\end{array}
\renewcommand{\arraystretch}{1}$$
(Here $T$ is the set of couples $(p, \lambda)$ where $1\leq p\leq r-1$ and $\lambda$ is a permutation for which $\lambda(1)<\cdots < \lambda (p)$ and $\lambda (p+1)<\cdots < \lambda (r)$; $(-1)^\lambda$ is the sign of $\lambda$.) 

A straightforward computation shows that the diagram
$$\xymatrix{
B\otimes UL\ar[d]^{d} \ar[rr]^\Delta && B\otimes (B\otimes UL)\ar[d]^{1\otimes d}\\
B\otimes UL \ar[rr]^{\Delta} && B\otimes (B\otimes UL)\,.}$$
commutes, and so 
  $(B\otimes UL,d)$ is a differential $B$-comodule.

From Fr\"oberg's formula it follows that $d: B\otimes (UL)(r)\to B\otimes (UL)(r+1)$. 
Since each $UL(p)$ is finite dimensional,     dual to the Fr\"oberg resolution is a resolution of $\mathbb Q$ by free $C$-modules:  
  $$\cdots \longrightarrow C\otimes \mbox{Hom}(UL(p+1),\mathbb Q)   \stackrel{d}{\longrightarrow}  
   C\otimes \mbox{Hom}(UL(p), \mathbb Q)\longrightarrow \cdots\,.$$
 This gives   isomorphisms  of graded vector spaces  Hom$((UL)(p), \mathbb Q))\cong \mbox{Tor}_p^C(\mathbb Q, \mathbb Q)$. 
 
 On the other hand, since $(\land V,d)$ is formal, it follows from (\cite{HaS}, \cite{S}) that $V$ is equipped with an extra gradation:
 \begin{enumerate}
 \item[$\bullet$] $V = \oplus_{p\geq 0} V(p)$,
 \item[$\bullet$] $d : V(0)\to 0$ and $d: V({p+1})\to (\land V)(p)$, where $(\land V)(p) = +_{p_1+\cdots + p_S = p} V({p_1}) \land \cdots \land V({p_s})$,
 \item[$\bullet$] The morphism $\pi : (\land V,d) \to H(\land V,d)$ given by $v\mapsto [v]$, $v\in V(0)$ and $V({\geq 1})\to 0$ is a quasi-isomorphism of bigraded spaces with $H(\land V,d) = H_0(\land V,d)$.
 \end{enumerate}
 
 Now the acyclic closure of $(\land V,d)$ can be constructed inductively by extending the acyclic closure of $(\land V({\leq p}),d)$ to one for $(\land V({\leq p+1}),d)$. This endows $U$ with a direct decomposition $U = \oplus_{p\geq 1} U(p)$, with
 $$d: U(p) \to (\land^+ V\otimes \land U)({p-1})\,,$$
 where $\land U$ and $\land V\otimes \land U$ have the obvious induced bigradations.
 
 In the quotient $H(\land V,d)\otimes_{(\land V,d)}\land U$ the differential has the form
 $$\cdots \to H(\land V,d) \otimes (\land U)({p+1})\to H(\land V,d)\otimes (\land U)(p)\to \cdots \to H(\land V,d)\otimes (\land U)(0)\to \mathbb Q\,.$$
 Since $\pi$ is a quasi-isomorphism so is $\pi\otimes id$, and so this is a free resolution of $\mathbb Q$ by $H(\land V,d)$-modules. In particular
 $$(\land U)(p)= \mbox{Tor}_{p}^{H(\land V,d)}(\mathbb Q,\mathbb Q)\,.$$
Since the graded algebras $C$ and $H(\land V,d)$ are isomorphic, the graded vector spaces $(\land U)({p})$ and $\mbox{Hom}\left( (UL)(p), \mathbb Q \right)$ are isomorphic. Since $(UL)= (UL)_0$ it follows that $U= U^0$ and $V= V^1$.  
 
 \hfill$\square$

 \section{The depth of a right-angled Artin group}
 
 The purpose of this section is the proof of Theorem 2 of the Introduction:

\vspace{3mm}\noindent {\bf Theorem 2.} {\sl Let $K$ be the classifying space for a   right-angled Artin group $A$ with generators $x_1, \dots , x_n$ and  with minimal Sullivan model $(\land V,d)$, and let
$$L = \mathbb L(x_i)_{1\leq i\leq n}/ I$$
be the corresponding weighted Lie algebra. Then
\begin{enumerate}
\item[(i)] The weighted Lie algebras $L_A$ and $L$ are isomorphic.
\item[(ii)] depth$\, (\land V,d)=$ depth$\, L=$ depth$\, UL=$ depth$\, \mathbb Q[A]\,.$
\item[(iii)] If $A$ is abelian then depth$\,\mathbb Q[A]= n$. Otherwise depth$\,\mathbb Q[A]=$ least $r+2$ for which there is a disconnecting simplex $\sigma$ with $\vert \sigma\vert = r$. 
\end{enumerate}}

\vspace{3mm} We first verify (i) and then through a sequence of results establish (ii).

\vspace{3mm}\noindent {\sl proof of Theorem 2(i).}  As observed in (\ref{UUn}), $L_A\cong E_V$. On the other hand, since $(\land V,d)$ satisfies (\ref{duo}) it follows that $V = \oplus_{p\geq 0} V(p)$ and that $V(0)\cong H^1(K)= C^1$.  
 
A basis of $ V(0)$ is given by the elements $t_1, \cdots , t_n$, and   we can take $x_1, \dots , x_n$ as the dual basis for $E_V(1)$, the dual of $sV(0)$. Thus   there is a 
natural surjective  Lie algebra morphism
 $$\pi : \mathbb L(x_1, \cdots , x_n) \to E_V$$
and by Proposition 3, $E_V$ is the quotient of $\mathbb L(x_1, \dots , x_n)$ by    the ideal generated by ker$\, \pi(2) : \mathbb L (x_1, \dots , x_n)(2)\to E_V(2)$.
 
 On the other hand, a basis of $V(1)$ is given by elements $t_{i,j}$, $i<j$, with $x_ix_j\neq x_jx_i$ and   $dt_{i,j} = t_i\land t_j$.   The corresponding elements $[x_i, x_j]\in \mathbb L(x_1, \dots , x_n)$   map to a basis of $E_V(2)$ while the other Lie brackets are in ker $\pi(2)$.  Therfore $\pi$ factors to give $L\stackrel{\cong}{\to} E_V$. 
 
  \hfill$\square$

 \vspace{3mm} \noindent {\sl proof of Theorem 2(ii).} Since $L\cong E_V$, Propositions 2 and 4  give
$$\mbox{depth}\, (\land V,d)= \mbox{depth}\, L=\mbox{depth}\, UL\,.$$
It remains to show that
\begin{eqnarray}\label{14}
\mbox{depth}\, UL = \mbox{depth}\, \mathbb Q[A]\,.
\end{eqnarray}
Before outlining and then detailing the proof of (\ref{14}) we   introduce some notation, establish a special case, and recall the fundamental result of Jensen and Meier \cite{JM}.

\vspace{3mm}
For any simplex $\sigma$ of an arbitrary polyhedron $R$, $\vert \sigma\vert$ denotes its dimension, $St\, \sigma := \mbox{star}\, \sigma$ is the set of simplices containing   $\sigma$, and  the closed star, $\overline{St\, \sigma}$ is the polyhedron consisting of the simplices in $St\, \sigma$ together with all their faces. The link, $Lk\, \sigma$ is the polyhedron of the simplices in $\overline{St\, \sigma}$ which contain no vertices of  $\sigma$.   The empty simplex $<\emptyset>$ has dimension $-1$ and $Lk <\emptyset>= R$. 

Finally, suppose $\sigma$ is a simplex in the flag polyhedron $P$. If $\tau= <y_1, \dots, y_s>$ is a simplex in $Lk\, \sigma$, where $\sigma = <x_1, \dots , x_r>$ then the vertices of $\tau$ commute with the vertices of $\sigma$ so that $\sigma*\tau = <x_1, \dots , x_r, y_1, \dots , y_s>$ is a simplex in $P$. In particular,
$$\overline{St\, \sigma} = \sigma * Lk\, \sigma\,.$$

Recall now that $x_1, \cdots , x_n$ are the generators of $A$. 
Suppose $\sigma = <y_1, \cdots , y_q>$ is a simplex in  $P$. If    $\tau = <z_1, \cdots , z_r>$ is a simplex in $Lk\, \sigma$ then the $z_i$ commute with the $y_j$ and so $\sigma*\tau = <y_1, \cdots , y_q, z_1, \cdots , z_r>$ is a simplex in $\overline{St\, \sigma}$. The vertices of $Lk\, \sigma$ and $\overline{St\, \sigma}$ then  generate sub right-angled Artin groups $A_{Lk\, \sigma}$ and $A_{\overline{St\, \sigma}}$  of $A$ whose flag polyhedra are then respectively $Lk\, \sigma$ and $\overline{St\,\sigma}$, and whose Lie algebras are respectively denoted $L_{Lk\, \sigma}$ and $L_{\overline{St\, \sigma}}$.

  \vspace{3mm}\noindent {\bf Lemma 3.}  Suppose  $\sigma = <x_{i_0}, \cdots , x_{i_r}>$ is a simplex in $P$. Then
 \begin{enumerate}
 \item[(i)] $$L_{\overline{St\, \sigma}}= \left(\bigoplus_{\lambda = 0}^r \mathbb Q x_{i_\lambda}\right) \, \times L_{Lk\, \sigma} \hspace{5mm}\mbox{and } A_{\overline{St\, \sigma}}=   \mbox{Ab}(  x_{i_\lambda})\times A_{Lk\, \sigma}$$
 where $\mbox{Ab} (x_{i_\lambda})$ is the free abelian group on $x_{i_\lambda}$.
 \item[(ii)] depth$\, UL_{\overline{St\, \sigma}}= (r+1)+ $ depth$\, UL_{Lk\, \sigma}$ and depth$\, \mathbb Q[A] = r+1 + $ depth$\, \mathbb Q[A_{Lk\, \sigma}]$
 \item[(iii)] depth$\, UL_{\overline{St\, \sigma}} =$ depth$\, \mathbb Q[A_{\overline{St\, \sigma}}]$ and depth$\, UL_{Lk\, \sigma}=$ depth$\, \mathbb Q[A_{Lk\, \sigma}]$.
 \end{enumerate}
 
 \vspace{3mm}\noindent {\sl proof.} (i) is essentially immediate from the definitions and (ii) follows at once.
 The second equality of (iii) follows by induction on the number of vertices,   and the first then follows from (ii). 

 \hfill$\square$

\vspace{3mm}\noindent {\bf Remark.} Lemma 3 establishes Theorem 2 if $P=\overline{St\, \sigma}$. In particular, if  $P$ is a simplex then $P = \overline{St\, x}$ for any vertex $x$, and so Theorem 2 is established.

\vspace{3mm}An identification of the depth of $\mathbb Q[A] $ is due to   Jensen and Meier \cite{JM}, and appears as a consequence of the following theorem. Denote by  the $\widetilde{H}(X;\mathbb Q)$ the reduced cohomology of $X$ with the convention that  $\widetilde{H}(\emptyset;\mathbb Q)= 0$.

\vspace{3mm}\noindent {\bf Theorem 5.} (\cite{JM}).  If the polyhedron $P$ is not a single simplex then
 $$\mbox{Ext}^p_{\mathbb Q[A]}(\mathbb Q, \mathbb Q[A]) = \bigoplus_{\sigma \subset P} \,\, \widetilde{H}^{p-\vert \sigma\vert -2}(Lk\, \sigma;\mathbb Q)\,,$$
 where the sum is over all the simplices $\sigma$ in $P$, including $\emptyset$.
 
 \vspace{3mm}\noindent {\bf Corollary 1.} {\sl If $P$ is not a single simplex, then $$\mbox{depth}\, \mathbb Q[A]= \mbox{greatest integer $m$ such that for each simplex $\sigma$ in $P$, }\widetilde{H}^{<m-\vert \sigma\vert -2}(\mbox{Lk}\, \sigma;\mathbb Q)= 0\,.$$}

\vspace{3mm}
Now, for any polyhedron $R$ we set
 $$n_R = \mbox{greatest integer $m$ such that for all simplices $\sigma$ in $R$}, \hspace{3mm} \widetilde{H}^{<m-\vert\sigma\vert -2}(Lk\, \sigma;\mathbb Q)= 0 \,.$$
 Thus the Corollary 1 can be restated as  
 $$\mbox{depth}\, \mathbb Q[A] = n_P\,,$$
if $P$ is not a simplex. Hence, 
   (\ref{14}) is equivalent to
 \begin{eqnarray}
 \label{prime}
 \mbox{depth}\, UL = n_P, \mbox{ if $P$ is not a simplex}.
 \end{eqnarray}

We first consider the case when $P$ is not connected:

\vspace{3mm}\noindent {\bf Lemma 4.}   {\sl If $P$ is not connected, then  depth$\, UL= n_P=1$.}

\vspace{3mm}\noindent {\sl proof.}  First of all, since $\mathbb Q[A]$ and $UL$ have no zero divisors, depth$\, \mathbb Q[A]\geq 1$, and depth$\, UL\geq 1$. Now since $P$ is not connected, $\widetilde{H}^0(\overline{St\, \emptyset}) \neq 0$, and so $n_P\leq 0-\vert \emptyset\vert + 2= 1$, so that $n_P= 1$.

On the other hand, since $P$ is not connected, $A$ is a non trivial free product $A = A'\# A''$. 
Let $(\land V',d)$ and $(\land V'',d)$ be the minimal Sullivan models for the classifying spaces of $A'$ and $A''$.  Then  (\cite[p.225]{FHTII}) the minimal Sullivan model, $  (\land V,d)$, for $P$ is the minimal Sullivan model of $(\land V',d)\oplus_{\mathbb Q} (\land V'',d)$, and \cite[Example 2 in section 10.2]{FHTII} gives   depth$\, UL=$ depth$(\land V,d) = 1$. \hfill$\square$

\vspace{3mm}
 It remains to prove (4) when $P$ is connected, which we do   by establishing separately the two inequalities depth$\, UL\leq n_P$ and depth$\, UL \geq n_P$.

 \vspace{3mm}\noindent {\bf Proof that depth$\, UL\leq n_P$}.

  We  show that if $\sigma$ is any $q$-simplex in $P$ then   $\widetilde{H}^{r}(Lk\,\sigma;\mathbb Q)$ embeds in Ext$_{UL}^{r+q+2}(\mathbb Q, UL)$. In particular,
$$\mbox{depth}\, UL\leq n_P\,.$$

For each simplicial set $X$ we denote by $C_*(X)$ and $C^*(X)$ the simplicial chain complex and the simplicial cochain complex on $X$ with rational coefficients. For recall $C_n(X)$ is the $\mathbb Q$-vector space generated by the $n$-simplices of $X$, and 
$$d<z_1, \dots, z_{n+1}> = \sum_{i=1}^{n+1}(-1)^{i+1} <z_1, \dots , \widehat{z_i},\dots , z_{n+1}>\,.$$
Then $C^*(X)= \mbox{Hom}(C_*(X), \mathbb Q)$ with the differential $(df)(\sigma) = (-1)^{deg\, f+1} f (d\sigma)$. 

 Now we   construct a chain map
$$\varphi_\sigma : \mbox{Hom}(C_*(Lk\, \sigma), \mathbb Q) \to \mbox{Hom}_{UL}(B_{*+\vert \sigma\vert +2}\otimes UL, UL)$$
and show that $H(\varphi_\sigma)$ induces an injection $\widetilde{H}^{r}(Lk\, \sigma;\mathbb Q) \to \mbox{Ext}^{r+q+2}_{UL}(\mathbb Q, UL)$.

When $\sigma = <\emptyset> $ we set
$$\varphi_{<\emptyset>} (f) <x_{i_1}\cdots   x_{i_{r+1}}>= f(<x_{i_1}  \cdots   x_{i_{r+1}}>)\cdot x_{i_1}\cdots x_{i_{r+1}}\,, \hspace{1cm} r\geq 0\,.$$
If $\sigma= <y_1, \cdots , y_{q+1}>$ we set
 $\varphi_\sigma (f)(\tau) = 0$ if $\tau \not\supset \sigma$ and, if $\tau \supset \sigma$, $$
\varphi_\sigma (f)(<x_{i_1}\cdots , x_{i_{r+1}}, y_1, \cdots , y_{q+1}> = f(<x_{i_1}, \cdots , x_{i_{r+1}}>)\cdot x_{i_1}\cdots x_{i_{r+1}}\,.$$

Note that $d\circ \varphi_\sigma (f)$ and $\varphi_\sigma (f\circ \partial)$ both vanish on simplices not containing $\sigma$ and, for the others, since the vertices of any complex commute in $UL$,  it follows from a simple computation that $(d\circ \varphi_\sigma)(f) =   \varphi_\sigma (f\circ \partial)$. Thus $\varphi_\sigma$ induces a linear map
$$H(\varphi_\sigma) = H^*(Lk\,\sigma;\mathbb Q) \to \mbox{Ext}_{UL}^{*+q+2}(\mathbb Q, UL)\,.$$

Now suppose that $f\in C^r(Lk(\sigma))$ is a cycle and $\varphi_\sigma (f)$ is a boundary. Then  there is a morphism $g : B_{r+q+1}\otimes UL\to UL$ such that $\varphi_\sigma (f) = g\circ d$.   By construction $\varphi_\sigma (f) (B_{r+q+2})\subset (UL)(r+1)$. Write $g = \sum_{i\geq 0} g_i$ with $g_i(B_{r+q+1})\subset (UL)(i)$. Now
$\varphi_\sigma (f) = g_{r}\circ d\,,$
and so we can   suppose that $g = g_{r}$.
Then, 
$$\renewcommand{\arraystretch}{1.7}
\begin{array}{ll}  f(<x_{i_1}, \cdots , x_{i_{r+1}}>)x_{i_1}\cdots x_{i_{r+1}} &=   d\circ g(<x_{i_1}, \cdots , x_{i_{r+1}}, y_1, \cdots , y_{q+1}>)\\
&= \displaystyle\sum_{j=1}^{r+1} (-1)^{j-1} g(<x_{i_1},\cdots \widehat{x_{i_j}}\cdots x_{i_{r+1}}, y_1, \cdots , y_{q+1}>)\cdot x_{i_j}\\
&+ \displaystyle\sum_{\ell = 1}^{q+1} (-1)^{\ell+r} g(<x_{i_1},\cdots   , x_{i_{r+1}}, y_1, \cdots \widehat{y_\ell} \cdots y_{q+1}>) \cdot y_\ell\,.
\end{array}
\renewcommand{\arraystretch}{1}$$

Write $$g<z_{i_1}, \dots , z_{i_{r}}, y_1, \ldots , y_{q+1}> = h(<z_{i_1}, \dots , z_{i_{r}}>)z_{i_1}\cdots z_{i_{r}} + \mu\,,$$
with $h(<z_{i_1}, \dots , z_{i_{r}}>)\in \mathbb Q$ and where $\mu$ is a linear combination of elements of $UL $ that either contains two times one of the variable $z_{i_j}$ or else contains a generator different of the $z_{i_j}$. 
Since there is no zero divisor in $UL$, this implies that
$$f(<x_{i_1}, \dots , x_{i_{r+1}}>) = \sum_{j=1}^{r+1} (-1)^{j-1} h(<x_{i_1}, \dots , \widehat{x_{i_j}}, \dots , x_{i_{r+1}}>)\,,$$
so $f$ is a boundary.
\hfill$\square$

\vspace{3mm}\noindent {\bf Proof that $n_P\leq \mbox{depth}\, UL$}.

We proceed by induction on the number of vertices in $P$. If $\sigma$ is a simplex in a sub polyhedron $R$ of $P$, we write $Lk (\sigma, R)$ for the link of $\sigma$ in $R$.

\vspace{3mm}\noindent {\bf Lemma 5.}  {\sl Suppose $\sigma$ is a simplex in $P$.
\begin{enumerate}
\item[(i)] If $\overline{St\, \sigma}$ is not a simplex then
$$\mbox{depth}\, UL_{\overline{St\, \sigma}}= n_{\overline{St\, \sigma}} \geq n_P\,.$$
\item[(ii)] If $\overline{St\, \sigma}$ is a simplex $\omega$ then for some simplex $\tau$, $\overline{St\, \tau}$ is not a simplex and 
$$\mbox{depth}\, UL_{\overline{St\, \sigma}} \geq \mbox{depth}\, UL_{\overline{St\, \tau}}\,.$$  
\end{enumerate}}

\vspace{3mm}\noindent {\sl proof.} 
(i) If $\sigma= \emptyset$ then $\overline{St\, \sigma} = P$ and $n_{\overline{St\, \sigma}} = n_P$. Otherwise,   $Lk\, \sigma$ is not empty or a single  simplex, because otherwise $\overline{St\, \sigma} = \sigma * (Lk\, \sigma)$ would be a simplex. Here we show that $$n_{\overline{St\, \sigma}} = n_{Lk\, \sigma} + \vert \sigma\vert +1\,.$$

Let $\tau\subset \overline{St\, \sigma}$ be a simplex. Then $\tau = \sigma_1*\tau_1$, with $\sigma_1\subset \sigma$, and $\tau_1\subset Lk\,\sigma$. Decompose $\sigma$ as $\sigma = \sigma_1*\sigma_2$. Then
$$Lk (\tau, \overline{St\, \sigma}) = \sigma_2*Lk(\tau_1, Lk\, \sigma)\,,$$
and so, $$\widetilde{H}^q(Lk(\tau, \overline{St\, \sigma})) = \widetilde{H}^{q-\vert \sigma_2\vert - 1}(Lk(\tau_1, Lk\, \sigma))\,.$$
Then
$$\begin{array}{ll}
n_{\overline{St\, \sigma}}&= \mbox{greatest integer $m$ such that for all simplex $\tau$, }\widetilde{H}^{<m-\vert \tau\vert -2}(Lk(\tau, \overline{St\, \sigma}))= 0\\&
=  \mbox{greatest integer $m$ such that  for all  $\tau$, }\widetilde{H}^{<m-\vert \tau\vert -\vert \sigma_2\vert - 3}(Lk(\tau_1, Lk\, \sigma))= 0\\&
=  \mbox{greatest integer $m$ such that  for all  $\tau$, }\widetilde{H}^{<m-\vert \sigma\vert -\vert\tau\vert + \vert \sigma_1\vert - 3}(Lk(\tau_1, Lk\, \sigma))= 0\\
&= \vert \sigma\vert + 1 +   \mbox{greatest integer $m$ such that  for all  $\tau_1$, }\widetilde{H}^{<m-\vert \tau_1\vert -2}(Lk(\tau_1, Lk\, \sigma))= 0\\
&= \vert \sigma\vert + 1 + n_{Lk\, \sigma}\,.\end{array}$$

Moreover, since $\overline{St\, \sigma}$ is not a simplex then $Lk\, \sigma$ is not a simplex. Thus by induction on the number of vertices, together with Lemma 3,
$$\mbox{depth}\, UL_{\overline{St\, \sigma}}= \vert \sigma\vert + 1 + \mbox{depth}\, UL_{Lk\, \sigma} = \vert \sigma\vert + 1+ n_{Lk\, \sigma} = n_{\overline{St\, \sigma}}\,.$$

On the other hand, for $\tau \subset Lk\, \sigma$, $Lk(\tau, Lk\, \sigma)= Lk(\tau*\sigma)$. Therefore $\widetilde{H}^k(Lk(\tau, Lk\, \sigma))= 0$ if $0<k<n_P-\vert \tau*\sigma\vert - 2$. This gives 
$$n_{Lk\, \sigma} \geq n_P-\vert \sigma\vert - 1$$
and so $$n_{\overline{St\, \sigma}} \geq n_P\,.$$

(ii)   If $\overline{St\, \sigma}$ is a simplex $\omega$ then $\overline{St\, \sigma} = \overline{St\, \omega}$ and so we may assume $\sigma = \overline{St}\, \sigma$,  in which case $\sigma$ is a maximal simplex. Write $\sigma= <x_1, \dots , x_r>$. For each vertex $y\not\in \sigma$, $(\overline{St\,y})\cap \sigma$ is a sub simplex $\sigma_y$ of $\sigma$. Since $\sigma$ is maximal, $\vert \sigma_y\vert <\vert \sigma\vert$. Now, choose a vertex $y$ with $\vert \sigma_y\vert$ maximal, and set $\tau = \sigma_y$. Then $Lk\, \tau$ is not connected, and so $\widetilde{H}^0(Lk\,\tau) \neq 0$. Therefore, 
$$n_{\overline{St\, \tau}} \leq 0 + \vert \tau\vert + 2 \leq \vert \sigma\vert + 1= \mbox{depth}\, UL_{\overline{St\, \sigma}}\,.$$
Thus by (i) depth$\, UL_{\overline{St\, \sigma}}\geq n_{\overline{St\, \tau}}= $ depth$\, UL_{\overline{St\, \tau}}$. 

\hfill$\square$

\vspace{3mm} Next, for each simplex $\sigma\subset P$ we may form the complex $B(\overline{St\, \sigma})\otimes UL_{\overline{St\, \sigma}}$. The inclusion $\overline{St\, \sigma}\subset P$ induces a morphism
$$\gamma : L_{\overline{St\, \sigma}} \to L\,,$$
and mapping the vertices $x_j\not\in \overline{St\, \sigma}$ to zero defines a retraction $L\to L_{\overline{St\, \sigma}}$. Thus $\gamma$ is an inclusion and $UL= UL_{\overline{St\, \sigma}}\otimes Z$ for some subspace $Z$. 

Consider the complex
$$B_*(\overline{St\, \sigma})\otimes UL : = \left( \, B_*(\overline{St\, \sigma})\otimes UL_{\overline{St\, \sigma}} \, \right) \otimes_{UL_{\overline{St\, \sigma}}} UL\,.$$
Then
$$H(\mbox{Hom}_{UL}(B( \overline{St\, \sigma})\otimes UL, UL) )= \mbox{Ext}_{UL_{\overline{St\, \sigma}}} (\mathbb Q, UL)$$
and it follows from Proposition 5 that 
\begin{eqnarray}
\label{star}
\mbox{depth}\, UL_{\overline{St\, \sigma}} = \mbox{least $p$ (or $\infty$) such that } H^p(\mbox{Hom}_{UL}(B(\overline{St\, \sigma})\otimes UL, UL)\neq 0\,.
\end{eqnarray}

\vspace{3mm}  \noindent {\bf Lemma 6.} {\sl If $P$ is connected, then for some simplex $\sigma\subset P$, depth$\, UL\geq $ depth$\, UL_{\overline{St\, \sigma}}$. }

\vspace{3mm}\noindent {\sl proof.} 
For any subpolyhedron $R\subset P$, let $W_*(R)$ denote the graded vector space in which    the $r$-simplices (including the empty simplex) in $R$ are a basis of $W_{r+1}(R)$.  We define a complex  
$$\dots C_p   \stackrel{d}{\longrightarrow} C_{p-1} \longrightarrow \dots C_1  \stackrel{\varepsilon}{\longrightarrow}  C_0 \to 0\,,$$ as follows:
First, set 
  $C_p = \oplus_{\vert \sigma\vert = p-1} W_*(\overline{St\, \sigma}) $.  Note that $C_0= W_*(\overline{St\, \emptyset})= W_*(P)$. Thus a basis of $C_p$ consists of the terms $(\sigma, \alpha)$ with $\vert\sigma\vert = p-1$ and $\alpha$  a simplex in $\overline{St\, \sigma}$, including the empty simplex. 

The differential $d$ is then defined by  
$$d(\sigma, \alpha) = \sum_{j=0}^p (-1)^j (<x_{i_0}, \dots , \widehat{x_{i_j}}, \dots , x_{i_p}>, \alpha)\,,$$ 
where  $\sigma= <x_{i_0}, \dots , x_{i_p}>$ and $\alpha \in \overline{St\, \sigma}$,   

This complex is exact. First of all $d : C_1\to C_0$ is surjective because  each simplex is in  $\overline{St\, x}$ for some $x$. 
Now let $\sum_{ij} \lambda_{ij} (\sigma_i, \alpha_j)$ be a cycle in some $C_p$, $\geq 1$. Then for each fixed $j$, $\sum_{i} \lambda_{ij}(\sigma_i, \alpha_{j})$ is a cycle. Since $\alpha_{j}\subset \overline{St\, \sigma_i}$ 
if and only if $\sigma_i \subset \overline{St\, \alpha_{j}}$, the complex generated by the $(\sigma, \alpha_{j})$   is  isomorphic to the usual chain complex of $\overline{St\, \alpha_{j}}$, and its homology   is zero, because the simplices $\sigma$ in  $\overline{St\, \alpha_j}$ include the empty simplex. This gives the acyclicity of the complex $(C_*,d)$. 

Next,  consider the double complex
$$\left(\, \oplus_{\vert\sigma\vert \geq 0}  \mbox{Hom}_{UL}(B_*(\overline{St\, \sigma})\otimes UL, UL), d_1+d_2\,\right)\,,$$
where $d_1$ is the internal differential in each   $\mbox{Hom}_{UL}(B_*(\overline{St\, \sigma})\otimes UL, UL)$. As a graded vector space  $ B_*(\overline{St\, \sigma}) = W_*(\overline{St\, \sigma})$ and $d_2= \mbox{Hom}_{UL}(d, -)$.
A standard computation shows that $d_1d_2+ d_2d_1= 0$. 
 
Consider first the differential $d_2$. The acyclicity of the complex $C_*$ shows that the natural morphism
$$(\mbox{Hom}_{UL}(B_*(P)\otimes UL, UL),d)\longrightarrow  \left(\, \oplus_{\vert\sigma\vert \geq 0}  \mbox{Hom}_{UL}(B_*(\overline{St\, \sigma})\otimes UL, UL), d_1+d_2\,\right) $$ is a quasi-isomorphism. In view of (\ref{star}) it follows that for some $\sigma$, depth$\, UL\geq $ depth$\, UL_{\overline{St\, \sigma}}$. 

\hfill$\square$

\vspace{3mm}

By Lemmas 5 and 6, for some simplex $\sigma$, $\overline{St\, \sigma}$ is not a simplex and 
$$\mbox{depth}\, UL\geq \mbox{depth}\, UL_{\overline{St\, \sigma}} = n_{\overline{St\, \sigma}} \geq n_P\,.$$
This proves depth$\, UL\geq n_P$ and there by Theorem 2(ii).
\hfill$\square$

 \vspace{3mm}\noindent {\sl proof of Theorem 2(iii).}  Lemma 6 shows that depth$\, UL\geq$ depth$\, UL_{\overline{St\, \sigma}}$ for some simplex $\sigma$, and by Lemma 5 we may suppose $\overline{St\, \sigma}$ is not a simplex. Since by Theorem 2(ii) we have depth$\, UL= n_P$ it follows from Lemma 5 that
$$\mbox{depth}\, UL= \mbox{depth}\, UL_{\overline{St\, \sigma}}\,.$$
Choose a $\sigma$ of maximal dimension with these properties.

Since $\overline{St\, \sigma}$ is not a simplex, $Lk\, \sigma$ is not empty and not a simplex. If $Lk\, \sigma$ is connected then
$$\mbox{depth}\, UL_{Lk\, \sigma} = \mbox{depth}\, UL_{\overline{St(\tau, Lk\, \sigma)}}$$
for some simplex $\tau\subset Lk\, \sigma$. But then
$$\overline{St\, \sigma} \supset \sigma *\tau * Lk(\tau, Lk\,\sigma) = \overline{St (\sigma*\tau)}$$
and depth$\, UL_{\overline{St(\sigma*\tau)}} =$ depth$\, UL_{\overline{St\, \sigma}}$, 
contrary to our hypothesis that dim$\, \sigma$ was maximal.

It follows that $Lk\, \sigma$ is not connected and thus
$$\mbox{depth}\, UL_{\overline{St\, \sigma}} = \vert \sigma\vert + 1+ \mbox{depth}\, UL_{Lk\, \sigma} = \vert \sigma\vert + 2\,.$$

On the other hand, for any simplex $\omega$, if $Lk\, \omega$ is not connected, then by Lemma 5,
$$\vert \omega\vert + 2= \mbox{depth}\, UL_{\overline{St\, \omega}} = n_{\overline{St\, \omega}} \geq n_P= \mbox{depth}\, UL\,.$$
\hfill$\square$

\vspace{3mm}\noindent {\bf Corollary.} {\sl A nilpotent right-angled Artin group $A$ is abelian.}

\vspace{3mm}\noindent {\sl proof.} Let $n$ be the number of generators of $A$, $L$ its Lie algebra and $(\land V,d)$ the minimal Sullivan model of its Salvetti complex $K$.  Since $A$ is nilpotent, dim$\,L<\infty$. By definition dim$\, L/[L,L]= n$. Therefore by \cite[Theorem 10.6]{FHTII},
$$\mbox{depth}\, L = \mbox{dim}\, L \,.$$

On the other hand, by \cite[Theorem 10.1]{FHTII} and \cite[Theorem 9.2]{FHTII},
$$\mbox{depth}\, (\land V,d) \leq \mbox{cat}\, (\land V,d) \leq \mbox{cat}\, K\,.$$

Now, since $K\subset (S^1)^n$, $K$ is a CW complex of dimension $\leq n$, and so cat$\, K\leq n$. All together this gives
	$$\mbox{dim}\, L/[L,L] = n \geq \mbox{cat}\, K \geq \mbox{depth}\, L= \mbox{dim}\, L\,.$$
 
Therefore $L$ is abelian.
\hfill$\square$

\vspace{3mm}\noindent {\bf Example 1.} Suppose $P$ is the union of simplices $\sigma$ and $\tau$ along a common proper simplex $\omega$. Then the vertices $x_0, \cdots , x_r$ of $\omega$ are central in $L$ and $L$ is the direct sum of the abelian Lie algebra they span together with the Lie algebra $F$ generated by the other vertices $y_1, \cdots , y_k$ of $\sigma$ and $z_1, \cdots , z_\ell$ of $\tau$. But $F$ is the free product of the Lie algebras $F'$ and $F''$ generated by the $y_i$ and the $z_j$.

Thus by Lemma 3,
$$\mbox{depth}\, UL =  \vert \omega\vert + 1 + \mbox{depth}\, UF\,.$$

Let  $A'$ and $A''$ be the right-angled Artin groups determined by the $y_i$ and the $z_j$. Then  $F$ is the weighted Lie algebra of $A'\#A''$, and by Lemma 4    depth$\, UF=$ depth$\, \mathbb Q[A'\#A''] = 1$. 
Therefore 
$$\mbox{depth}\, UL= \vert \omega\vert + 2\,.$$

\vspace{3mm}\noindent {\bf Example 2.} Suppose $A$ has generators $x_i, y_j$ where the relations consist of relations among the $x_i$, relations among the $y_j$ and one additional relation $x_1y_1= y_1x_1$. If there are at least 2 of the $x_i$ and of the $y_j$ and if the graph of $A$ is connected then depth$\, A \leq 3  $. 

 In fact, let $\sigma = <x_1, y_1>$. Then $Lk\, \sigma$ is disconnected and so $\widetilde{H}^0(Lk\, \sigma)\neq 0$. Thus it follows from the formula of Jensen and Meier that $n_P-\vert \sigma\vert -2\leq 0$ and so depth$\, K= n_P\leq 3$.

 \vspace{3mm}
 Riemann surfaces of genus $g\geq 1$, nilmanifolds and classifying spaces of right-angled Artin groups are all $K(G,1)$ spaces whose minimal model $(\land V,d)$ satisfies $V = V^1$. This suggests the following question.

 \vspace{2mm}\noindent {\bf Problem.} Suppose that a finite CW complex $X$ is a $K(G,1)$ and that the minimal Sullivan model of $X$, $(\land V,d)$ satisfies $V = V^1$. Denote by $L$ the Lie algebra associated to the rational MalcevLie algebra $L_G$. Is there more generally a relation between the depth of $UL$ and   of $\mathbb Q[G]$ ?

   \vspace{1cm} Institut de Recherche en Math\'ematique et en Physique, Universit\'e Catholique de Louvain, 2, Chemin du Cyclotron, 1348 Louvain-La-Neuve, Belgium, yves.felix@uclouvain.be
   
   \vspace{3mm} Department of Mathematics, Mathematics Building, University of Maryland, College Park, MD 20742, United States, shalper@umd.edu
   
\end{document}